\documentclass[11pt]{amsart}
\usepackage[latin1]{inputenc}

\usepackage[mathcal]{eucal}
\usepackage{amsmath}
\usepackage{amsfonts}
\usepackage{amssymb}
\usepackage{amsthm}
\usepackage{graphicx,epsfig}
\usepackage{amscd}
\DeclareMathOperator{\argsinh}{argsinh}

\DeclareMathOperator{\re}{Re}
\DeclareMathOperator{\im}{Im}

\newcommand{\nil}{\mathrm{Nil}_3}

\newcommand{\D}{\mathbb{D}}

\newcommand{\R}{\mathbb{R}}
\newcommand{\h}{\mathbb{H}}

\newcommand{\C}{\mathbb{C}}
\newcommand{\s}{\mathbb{S}}
\newcommand{\LL}{\mathbb{L}}

\newcommand{\rmd}{\mathrm{d}}

\newcommand{\cS}{{\mathcal S}}

\newcommand{\cP}{{\mathcal P}}

\begin{document}

\newtheorem{thm}{Theorem}[section]
\newtheorem*{thmintro}{Theorem}
\newtheorem{cor}[thm]{Corollary}
\newtheorem{prop}[thm]{Proposition}
\newtheorem{app}[thm]{Application}
\newtheorem{lemma}[thm]{Lemma}
\newtheorem{notation}[thm]{Notations}
\newtheorem{hypothesis}[thm]{Hypothesis}

\newtheorem{defin}[thm]{Definition}
\newenvironment{defn}{\begin{defin} \rm}{\end{defin}}
\newtheorem{remk}[thm]{Remark}
\newenvironment{rem}{\begin{remk} \rm}{\end{remk}}
\newtheorem{exa}[thm]{Example}
\newenvironment{ex}{\begin{exa} \rm}{\end{exa}}

\title[The Gauss map of minimal surfaces]{The Gauss map of minimal surfaces \\
in the Heisenberg group}
\author{Beno\^\i t Daniel}
\date{}


\address{Universit\'e Paris-Est Cr\'eteil, D\'epartement de Math\'ematiques, UFR des Sciences et Technologies, 61 avenue du G\'en\'eral de Gaulle, 94010 Cr\'eteil cedex, FRANCE}
\email{daniel@univ-paris12.fr}

\begin{abstract}
We study the Gauss map of minimal surfaces in the Heisenberg group $\nil$ endowed with a left-invariant Riemannian metric.
We prove that the Gauss map of a nowhere vertical minimal surface is harmonic into the hyperbolic plane $\h^2$. Conversely, any nowhere antiholomorphic harmonic map into $\h^2$ is the Gauss map of a nowhere vertical minimal surface. Finally, we study the image of the Gauss map of complete nowhere vertical minimal surfaces.
\end{abstract}

\maketitle

\section{Introduction}

The Gauss map of minimal and constant mean curvature (CMC) surfaces has been the object of a huge amount of investigations. For example, the Gauss map of minimal surfaces in Euclidean space $\R^3$ and of CMC $1$ surfaces in hyperbolic space $\h^3$ is meromorphic, the Gauss map of CMC surfaces in $\R^3$ (respectively, spacelike CMC in Minkowski space $\LL^3$) is harmonic into the sphere $\s^2$ (respectively, the hyperbolic plane $\h^2$). Hence a lot of results on the geometry of minimal and CMC surfaces have been obtained using the theory of meromorphic and harmonic maps.

More recently, I. Fern\'andez and P. Mira proved the existence of a harmonic ``hyperbolic Gauss map" into $\h^2$ for CMC $\frac12$ surfaces in $\h^2\times\R$ that are nowhere vertical (\cite{fernandezmira}). They also proved that any harmonic map (under a hypothesis) can be the hyperbolic Gauss map of a CMC $\frac12$ immersion.

Because of the local isometric correspondence between CMC $\frac12$ surfaces in $\h^2\times\R$ and minimal surfaces in the $3$-dimensional Heisenberg group $\nil$ (see \cite{sister}), it is natural to study the same problem for minimal surfaces in $\nil$ (endowed with a left-invariant Riemannian metric). The Lie group $\nil$ is a $3$-dimensional homogeneous manifold with a $4$-dimensional isometry group; hence it is one of the most simple $3$-manifolds apart from space-forms. Moreover, it is a Riemannian fibration over the Euclidean plane $\R^2$.

Identifying (by left multiplication) the tangent space of $\nil$ at any point with the Lie algebra of $\nil$, the Gauss map of a surface can be considered as a map into the unit sphere of the Lie algebra. If the surface is nowhere vertical, i.e., nowhere tangent to the fibers, then the image of the Gauss map is contained in the northern hemisphere (up to a change of orientation).
We first prove that the Gauss map of a nowhere vertical minimal surface in $\nil$ is harmonic into the northern hemisphere endowed with the \emph{hyperbolic} metric (i.e., of constant curvature $-1$; theorem \ref{gharmonic} and corollary \ref{corgharmonic}). 

Conversely, we prove that any nowhere antiholomorphic harmonic map into $\h^2$ is the Gauss map of a nowhere vertical minimal surface in $\nil$ (theorem \ref{weierstrass}). For this purpose we obtain a Weierstrass-type representation in terms of the harmonic Gauss map. We observe that the formulas turn out to be much simpler than the corresponding ones for CMC $\frac12$ surfaces in $\h^2\times\R$ (\cite{fernandezmira}).

In section \ref{sectioncomplete}, we get some properties of the image of the Gauss of complete nowhere vertical minimal surfaces in $\nil$, using the theory of harmonic maps into $\h^2$ and of CMC surfaces in $\LL^3$ (theorem \ref{gaussianimage}). Finally, in section \ref{sectionexamples} we give the Gauss map of the classical examples of minimal surfaces in $\nil$ and we construct a new (and explicit) example of an entire minimal graph in $\nil$ (example \ref{semitrough}).

In the present paper, the model used for the hyperbolic plane is the Poincar\'e disk, i.e., $$\h^2=\D=\{z\in\C;|z|<1\}$$ endowed with the metric
$$\frac{4|\rmd z|^2}{(1-|z|^2)^2}.$$
Also, we denote by $\bar\C=\C\cup\{\infty\}$ the Riemann sphere.

\section{The Heisenberg group $\nil$} \label{heisenberg}

The $3$-dimensional Heisenberg group $\nil$ can be viewed as $\R^3$ endowed with the metric
$$\rmd x_1^2+\rmd x_2^2+\left(\frac12(x_2\rmd x_1-x_1\rmd x_2)+\rmd x_3\right)^2.$$
The projection $\pi:\nil\to\R^2,(x_1,x_2,x_3)\mapsto(x_1,x_2)$ is a Riemannian fibration; we will identify $\R^2$ and $\C$.

We consider the left-invariant orthonormal frame $(E_1,E_2,E_3)$ defined by
$$E_1=\frac\partial{\partial x_1}-\frac{x_2}2\frac\partial{\partial x_3},\quad
E_2=\frac\partial{\partial x_2}+\frac{x_1}2\frac\partial{\partial x_3},\quad
E_3=\frac\partial{\partial x_3}.$$ We call it the canonical frame.

The expression of the Riemannian connection $\hat\nabla$ in this frame is the following:
$$\hat\nabla_{E_1}E_1=0,\quad
\hat\nabla_{E_2}E_1=-\frac12E_3,\quad
\hat\nabla_{E_3}E_1=-\frac12E_2,$$
$$\hat\nabla_{E_1}E_2=\frac12E_3,\quad
\hat\nabla_{E_2}E_2=0,\quad
\hat\nabla_{E_3}E_2=\frac12E_1,$$
$$\hat\nabla_{E_1}E_3=-\frac12E_2,\quad
\hat\nabla_{E_2}E_3=\frac12E_1,\quad
\hat\nabla_{E_3}E_3=0.$$

A vector is said to be vertical if it is proportional to $E_3$, and horizontal if it is orthogonal to $E_3$. A surface is said to be nowhere vertical if $E_3$ is nowhere tangent to it, i.e., if the restriction of $\pi$ to the surface is a local diffeomorphism.

We identify the tangent space at a point $x=(x_1,x_2,x_3)$ with the Lie algebra of $\nil$ by left multiplication by $x^{-1}$; in other words, we can identify a vector at $x$ with its coordinates in the frame $(E_1,E_2,E_3)$. These coordinates will be denoted into brackets; we then have
$$a_1\frac\partial{\partial x_1}+a_2\frac\partial{\partial x_2}+a_3\frac\partial{\partial x_3}=
\left[\begin{array}{c}
a_1 \\
a_2 \\
a_3+\frac12(x_2a_1-x_1a_2)
      \end{array}\right].$$

The isometry group of $\nil$ is $4$-dimensional and has two connected components: isometries preserving orientation of the fibers and the base of the fibration, and those reversing both of them. 
The translations along the $x_1$-axis are given by $$(x_1,x_2,x_3)\mapsto\left(x_1+t,x_2,x_3+\frac{tx_2}2\right)$$ with $t\in\R$, those along the $x_2$-axis by $$(x_1,x_2,x_3)\mapsto\left(x_1,x_2+t,x_3-\frac{tx_1}2\right),$$ and those along the $x_3$-axis by $$(x_1,x_2,x_3)\mapsto(x_1,x_2,x_3+t).$$ The canonical frame is invariant by translations. A rotation by $\theta$ about a vertical fiber in $\nil$ induces in the Lie algebra the rotation of angle $\theta$ about $E_3$. For more details, see \cite{mercuri}. 

\section{The Gauss map} \label{sectiongauss}

Let $\Sigma$ be an oriented Riemann surface and $z=u+iv$ a conformal coordinate in $\Sigma$. Let $X:\Sigma\to\nil$ be a conformal immersion. We denote by $F=\pi\circ X$  the horizontal projection of $X$ and by $h:\Sigma\to\R$ the third coordinate of $X$ in the model described in section \ref{heisenberg}; we have $X=(F,h)$ in this model. We regard $F$ as a complex-valued function, identifying $\C$ with $\R^2$. We denote by $N:\Sigma\to\s^2$ the unit normal to $X$, where $\s^2$ is the unit sphere of the Lie algebra of $\nil$.

\begin{defn}
The Gauss map of $X$ is the map $g=\varphi\circ N:\Sigma\to\bar\C$ where $\varphi$ is the stereographic projection with respect to the southern pole, i.e., $g:\Sigma\to\bar\C$ is defined by
$$N=\frac1{1+|g|^2}\left[\begin{array}{c} 2\re g \\ 2\im g \\ 1-|g|^2
          \end{array}\right].$$
\end{defn}

\begin{defn}
We set
$$\eta=2\langle E_3,X_z\rangle.$$ We call the pair $(g,\eta)$ the Weierstrass data of the immersion $X$.
\end{defn}

\begin{rem}
The other choice of normal would replace $g$ by $\tilde g=-\frac1{\bar g}$.
\end{rem}

We have $$X_z=\frac12\left[\begin{array}{c} (F+\bar F)_z \\ i(\bar F-F)_z \\ \eta \end{array}\right],\quad
X_{\bar z}=\frac12\left[\begin{array}{c} (F+\bar F)_{\bar z} \\ i(\bar F-F)_{\bar z} \\ \bar\eta \end{array}\right].$$

Observe that the conformality of $X$ is equivalent to 
\begin{equation} \label{conformality}
F_z\bar F_z=-\frac{\eta^2}4.
\end{equation}

We also have $$X_z\times X_{\bar z}=\frac i2\left[\begin{array}{c}
\re(\eta F_{\bar z}-\bar\eta F_z) \\
\im(\eta F_{\bar z}-\bar\eta F_z) \\
|F_z|^2-|F_{\bar z}|^2
                                                          \end{array}\right],$$
and thus
$$X_u\times X_v=-2iX_z\times X_{\bar z}=\left[\begin{array}{c}
\re(\eta F_{\bar z}-\bar\eta F_z) \\
\im(\eta F_{\bar z}-\bar\eta F_z) \\
|F_z|^2-|F_{\bar z}|^2
                                                          \end{array}\right].$$

We first notice that
\begin{itemize}
\item $g=0$ or $\infty$ if and only if $\eta=0$,
\item $g=0$ if and only if $F_{\bar z}=0$,
\item $g=\infty$ if and only if $F_z=0$.
\end{itemize}

We compute that
$$||X_u\times X_v||=2\langle X_z,X_{\bar z}\rangle=|F_z|^2+|F_{\bar z}|^2+\frac12|\eta|^2.$$
From this we get, when $g\neq\infty$,
$$\frac{2g}{1+|g|^2}=\frac{\eta F_{\bar z}-\bar\eta F_z}{|F_z|^2+|F_{\bar z}|^2+\frac12|\eta|^2},\quad
\frac{1-|g|^2}{1+|g|^2}=\frac{|F_z|^2-|F_{\bar z}|^2}{|F_z|^2+|F_{\bar z}|^2+\frac12|\eta|^2},$$
and so
$$g=\frac{\eta F_{\bar z}-\bar\eta F_z}{2|F_z|^2+\frac12|\eta|^2}.$$
Using \eqref{conformality}, we conclude that, when $g\neq\infty$,
\begin{equation} \label{diffF}
\bar gF_z=-\frac\eta2,\quad F_{\bar z}=\frac{g\bar\eta}2.
\end{equation}

\begin{rem} \label{extension}
This shows that the functions $\frac{\eta}{\bar g}$ and $g\bar\eta$ can be extended smoothly to points where $g=0$ or $\infty$.
\end{rem}

Formulas \eqref{conformality} and \eqref{diffF} characterize $(g,\eta)$ up to the transformation $(g,\eta)\mapsto(-g,-\eta)$.
Observe that $F$, $g$ and $\eta$ are defined independently of the model of $\nil$, whereas $h$ depends on the model. This function $h$ satisfies
\begin{equation} \label{diffh}
h_z=\frac12\eta-\frac i4(\bar FF_z-F\bar F_z).
\end{equation}

We first notice the following fact.

\begin{prop} \label{gnotzero}
Let $X:\Sigma\to\nil$ be a conformal immersion. Then its Gauss map $g$ cannot be identically equal to $0$ or $\infty$ on an open set. 
\end{prop}

\begin{proof}
Assume that $g\equiv 0$ on an open set $U$. Then by \eqref{diffF} we have $\eta\equiv 0$ and $F_{\bar z}\equiv 0$ on $U$, and thus $F_{z\bar z}\equiv 0$ on $U$. Hence by \eqref{diffh} we get $h_z=-\frac i4\bar FF_z$, and so $h_{z\bar z}=-\frac i4|F_z|^2$. Since $h_{z\bar z}\in\R$ we conclude that $F_z\equiv 0$ on $U$, and so $X_z\times X_{\bar z}\equiv 0$, which is impossible, since $X$ is an immersion.

In the same way, $g$ cannot be identically equal to $\infty$ on an open set. 
\end{proof}

\begin{thm}  \label{gharmonic}
Let $X:\Sigma\to\nil$ be a conformal minimal immersion. Then its Gauss map $g:\Sigma\to\bar\C$ satisfies
\begin{equation} \label{eqg}
(1-|g|^2)g_{z\bar z}+2\bar gg_zg_{\bar z}=0.
\end{equation}
This equality holds when $g\neq\infty$; in the neighbourhood of a point where $g=\infty$, $g$ needs to be replaced by $\frac1g$.
\end{thm}

\begin{cor} \label{corgharmonic}
Let $X:\Sigma\to\nil$ be a conformal minimal immersion that is nowhere vertical. Assume that its normal is upward. Then its Gauss map $g$ takes values in the hyperbolic disk $\h^2$, and it is harmonic (for the hyperbolic metric).
\end{cor}

\begin{proof}[Proof of corollary \ref{corgharmonic}]
The third component of $N$ is $\frac{1-|g|^2}{1+|g|^2}$, so $N$ is upward if and only if $|g|<1$, i.e., $g\in\h^2$. Also, the equation
$$g_{z\bar z}+\frac{2\bar g}{1-|g|^2}g_zg_{\bar z}=0$$ means that $g$ is a harmonic map into $\h^2$.
\end{proof}

Before proving the theorem, we translate the minimality condition in terms of $F$, $g$ and $\eta$.

\begin{lemma} \label{minimalitycondition}
The conformal immersion $X$ is minimal if and only if
\begin{equation} \label{systemminimal}
\left\{\begin{array}{c}
F_{z\bar z}=\frac i8|\eta|^2\frac{|g|^2-1}{|g|^2}g, \\
\eta_{\bar z}+\bar\eta_z=0.
         \end{array}\right.
\end{equation}
The first equation in \eqref{systemminimal} holds when $g\neq0,\infty$. When $g=0,\infty$ it means $F_{z\bar z}=0$ (see also remark \ref{extension}).
\end{lemma}

\begin{proof}
The conformal immersion $X$ is minimal if and only if $$\hat\nabla_{X_u}X_u+\hat\nabla_{X_v}X_v=0.$$ We have
$$\hat\nabla_{X_u}X_u=\frac14\left[\begin{array}{c}
2(F+\bar F)_{uu}+i(\eta+\bar\eta)(\bar F-F)_u \\
2i(\bar F-F)_{uu}-(\eta+\bar\eta)(F+\bar F)_u \\
2(\eta+\bar\eta)_u
                            \end{array}\right],$$
$$\hat\nabla_{X_v}X_v=\frac14\left[\begin{array}{c}
2(F+\bar F)_{vv}-(\eta-\bar\eta)(\bar F-F)_v \\
2i(\bar F-F)_{vv}-i(\eta-\bar\eta)(F+\bar F)_v \\
2i(\eta-\bar\eta)_v
                            \end{array}\right].$$
Hence the horizontal part of $\hat\nabla_{X_u}X_u+\hat\nabla_{X_v}X_v$ vanishes if and only if
$$F_{z\bar z}=\frac i4(\bar\eta F_z+\eta F_{\bar z}),$$
and its vertical part vanishes if and only if
$$\eta_{\bar z}+\bar\eta_z=0.$$
When $g\neq0,\infty$, we conclude using \eqref{diffF}. When $g=0,\infty$, we have $\eta=0$ and so $F_{z\bar z}=0$.
\end{proof}

\begin{proof}[Proof of theorem \ref{gharmonic}]
We first restrict ourselves to a domain on which $g\neq\infty$, i.e., on which $F_z\neq 0$. Using \eqref{diffF} we have $g^2=-\frac{F_{\bar z}}{\bar F_{\bar z}}$, and thus
$$2gg_z=-\frac{F_{z\bar z}}{\bar F_{\bar z}}+\frac{F_{\bar z}\bar F_{z\bar z}}{\bar F_{\bar z}^2},$$ and by \eqref{systemminimal} we obtain, when $g\neq 0$,
\begin{equation} \label{gz}
g_z=\frac i4F_z(1-|g|^2)^2.
\end{equation}
By continuity this equality also holds when $g=0$.
Differentiating \eqref{gz} with respect to $\bar z$, we conclude, using \eqref{diffF}, \eqref{systemminimal} and \eqref{gz}, that $g$ satisfies \eqref{eqg}.

In the same way, in the neighbourhood of a point where $g=\infty$, we prove that $\frac1g$ satisfies \eqref{eqg}.
\end{proof}

\begin{rem} \label{verticalplane}
If $|g|\equiv1$ on an open set, then by \eqref{eqg} we obtain that $g$ is constant, and so, using \eqref{diffF}, we obtain that $\bar gF+g\bar F$ is constant. Hence $X(\Sigma)$ is a vertical plane, i.e., a surface of equation $ax_1+bx_2=c$ for some $(a,b)\neq(0,0)$ and some $c\in\R$ (such a surface is minimal and flat).
\end{rem}

\begin{prop}
Let $X:\Sigma\to\nil$ be a conformal minimal immersion such that $X(\Sigma)$ is not a vertical plane. Let $(g,\eta)$ be its Weierstrass data. Then we have
\begin{equation} \label{eta}
\eta=8i\frac{\bar gg_z}{(1-|g|^2)^2},
\end{equation}
and the metric induced on $\Sigma$ by $X$ is
\begin{equation} \label{metric}
\rmd s^2=16\frac{(1+|g|^2)^2}{(1-|g|^2)^4}|g_z|^2|\rmd z|^2.
\end{equation}
These expressions hold at points where $|g|\neq 1,\infty$, and extend smoothly at points where $|g|=1,\infty$. In particular, if $|g|<1$, then $g$ is nowhere antiholomorphic (i.e., $g_z$ does not vanish).
\end{prop}

\begin{proof}
Formula \eqref{eta} comes from \eqref{diffF} and \eqref{gz}.

Uusing \eqref{diffF} and \eqref{eta}, we get
$$\langle X_z,X_{\bar z}\rangle=\frac12(|F_z|^2+|F_{\bar z}|^2)+\frac14|\eta|^2
=8\frac{(1+|g|^2)^2}{(1-|g|^2)^4}|g_z|^2|\rmd z|^2,$$
so we obtain \eqref{metric}.

Finally, we have $\langle X_z,X_{\bar z}\rangle>0$, since $X$ is an immersion; hence, if $|g|<1$, this shows that $g$ is nowhere antiholomorphic.
\end{proof}

\begin{rem}
A. Sanini proved that the Gauss map, regarded as a map into the Grassmann bundle of the $2$-planes of the tangent bundle of $\nil$, of a surface in $\nil$ is conformal if and only if the surface is minimal (\cite{sanini}).
\end{rem}

\section{Minimal immersions with prescribed Gauss map} \label{sectionprescribed}

\begin{thm} \label{weierstrass}
Let $\Sigma$ be a simply-connected Riemann surface.
Let $g:\Sigma\to\h^2$ be a harmonic map that is nowhere antiholomorphic. Let $z_0\in\Sigma$, $F_0\in\C$ and $h_0\in\R$.

Then there exists a unique conformal minimal immersion $X:\Sigma\to\nil$ such that $g$ is the Gauss map of $X$ and $X(z_0)=(F_0,h_0)$.

Moreover, the immersion $X=(F,h)$ satisfies
$$F_z=-4i\frac{g_z}{(1-|g|^2)^2},\quad
F_{\bar z}=-4i\frac{g^2\bar g_{\bar z}}{(1-|g|^2)^2},$$
$$h_z=4i\frac{\bar gg_z}{(1-|g|^2)^2}-\frac i4(\bar FF_z-F\bar F_z).$$
\end{thm}

\begin{proof}
We first recover the horizontal part $F$ using the following fact. We claim that
the differential system
\begin{equation} \label{systemF}
\left\{\begin{array}{ccc}
F_z & = & -4i\frac{g_z}{(1-|g|^2)^2} \\
F_{\bar z} & = & -4i\frac{g^2\bar g_{\bar z}}{(1-|g|^2)^2}
      \end{array}\right.
\end{equation}
has a unique solution $F:\Sigma\to\C$ satisfying $F(z_0)=F_0$. Indeed, setting
$A=-4i\frac{g_z}{(1-|g|^2)^2}$ and $B=-4i\frac{g^2\bar g_{\bar z}}{(1-|g|^2)^2}$, it suffices to check that $A_{\bar z}=B_z$ (since $\Sigma$ is simply-connected). Using the fact that $g$ is harmonic, we get that $A_{\bar z}=B_z=-8i\frac{gg_z\bar g_{\bar z}}{(1-|g|^2)^3}$, which proves the claim.

We now have to recover the function $h$. We set 
\begin{equation} \label{definitioneta}
\eta=8i\frac{\bar gg_z}{(1-|g|^2)^2}.
\end{equation}
We claim that the differential equation
\begin{equation} \label{eqh}
h_z=\frac12\eta-\frac i4(\bar FF_z-F\bar F_z)
\end{equation}
has a unique solution $h:\Omega\to\R$ satisfying $h(z_0)=h_0$. Indeed, it suffices to check that
$$\left(\frac12\eta-\frac i4(\bar FF_z-F\bar F_z)\right)_{\bar z}\in\R.$$
Using the fact that $g$ is harmonic and \eqref{systemF}, we get
$$\left(\frac12\eta-\frac i4(\bar FF_z-F\bar F_z)\right)_{\bar z}=\frac i4(F\bar F_{z\bar z}-\bar FF_{z\bar z})\in\R,$$
which proves the claim.

We now set $X=(F,h):\Sigma\to\nil$, and check that $X$ satisfies the conclusions of the theorem. It is clear that $X(z_0)=(F_0,h_0)$.

By \eqref{eqh}, we have $\eta=2\langle E_3,X_z\rangle$; moreover, by \eqref{systemF} and \eqref{definitioneta}, \eqref{conformality} is satisfied, so $X$ is a conformal map. Also, $\langle X_z,X_{\bar z}\rangle=\frac12(|F_z|^2+|F_{\bar z}|^2)+\frac14|\eta|^2>0$, since $g_z$ does not vanish; hence $X$ is an immersion. We also observe that \eqref{diffF} holds, which means that $g$ is the Gauss map of $X$. Finally, \eqref{systemminimal} holds, so $X$ is minimal.
\end{proof}

\begin{rem}
If we start with a map $g:\Sigma\to\bar\C$ satisfying equation \eqref{eqg} (without assuming that $|g|<1$), then it may be impossible to recover a minimal immersion $X:\Sigma\to\nil$. However, in some cases it is possible.
\begin{itemize}
\item If $|g|\equiv1$ on an open set, then we recover a vertical plane (remark \ref{verticalplane}).
\item Assume that $|g|$ is not identically $1$ on any open set. Let $\Sigma_0$ be a connected domain in $\Sigma$ on which the functions $-4i\frac{g_z}{(1-|g|^2)^2}$ and $-4i\frac{\bar g^2g_z}{(1-|g|^2)^2}$ can be extended smoothly to a pair of functions $A$ and $B$ that do not vanish simultaneously. Let $\tilde\Sigma_0$ be the universal cover of $\Sigma_0$. Then we can recover an immersion $X:\tilde\Sigma_0\to\nil$ (the proof is the same as that for theorem \ref{weierstrass}; the fact that $A$ and $B$ do not vanish simultaneously ensures that $X$ is an immersion). See the examples in section \ref{sectionexamples}.
\end{itemize}
\end{rem}

\begin{rem}
Such representation theorems involving harmonic maps were proved by K. Kenmotsu for CMC surfaces in $\R^3$ (\cite{kenmotsu}), by K. Akutagawa and S. Nishikawa for spacelike CMC surfaces in $\LL^3$ (\cite{akni}) and by I. Fern\'andez and P. Mira for CMC $\frac12$ surfaces in $\h^2\times\R$ (\cite{fernandezmira}).
\end{rem}

\begin{rem}
F. Mercuri, S. Montaldo and P. Piu also obtained a Weier-strass-type representation for minimal surfaces in $\nil$ (\cite{mercuri2}). However, their representation involves a pair of functions satisfying a rather complicated system of PDE, whose solutions are difficult to find. On the contrary, theorem \ref{weierstrass} allows to construct a minimal surface starting simply with a harmonic map into $\h^2$. Examples are provided in section \ref{sectionexamples}.
\end{rem}

\section{The Hopf differential}

The Hopf differential of a harmonic map $g:\Sigma\to\h^2$ is
$$Q=\frac4{(1-|g|^2)^2}g_z\bar g_z\rmd z^2.$$ It is well known that it is holomorphic. 

On the other hand, U. Abresch and H. Rosenberg proved the existence of a holomorphic quadratic differential for CMC surfaces in $\nil$, and more generally in other homogeneous $3$-manifolds (\cite{ar1}, \cite{ar2}); this quadratic differential is a complex linear combination of the $(2,0)$ part of the second fundamental form of the surface and of $\eta^2\rmd z^2$.

\begin{prop} \label{hopf}
The Abresch-Rosenberg differential of a nowhere vertical minimal surface in $\nil$ coincides (up to a constant) with the Hopf differential of its Gauss map.
\end{prop}

\begin{proof}
We use the same notation as in sections \ref{sectiongauss} and \ref{sectionprescribed}.

We have $$\hat\nabla_{X_u}X_v=\frac14\left[\begin{array}{c}
2(F+\bar F)_{uv}+\frac i2(\eta+\bar\eta)(\bar F-F)_v-\frac12(\eta-\bar\eta)(\bar F-F)_u \\
2i(\bar F-F)_{uv}-\frac12(\eta+\bar\eta)(F+\bar F)_v-\frac i2(\eta-\bar\eta)(F+\bar F)_u \\
i(\eta-\bar\eta)_u+(\eta+\bar\eta)_v
                            \end{array}\right].$$
The last coordinate is obtained, using that it is equal to $2i(\eta-\bar\eta)_u+i(F_u\bar F_v-\bar F_uF_v)$ and to $2(\eta+\bar\eta)_v-i(F_u\bar F_v-\bar F_uF_v)$.
Using a computation done in the proof of lemma \ref{minimalitycondition}, we get
\begin{eqnarray*}
\lefteqn{\hat\nabla_{X_u}X_u-\hat\nabla_{X_v}X_v-2i\hat\nabla_{X_u}X_v} \\
& = & \left[\begin{array}{c}
2(F+\bar F)_{zz}+\frac i2\eta(\bar F-F)_u+\frac12\eta(\bar F-F)_v \\
2i(\bar F-F)_{zz}+\frac i2\eta(F+\bar F)_v-\frac12\eta(F+\bar F)_u \\
2\eta_z
                            \end{array}\right].\end{eqnarray*}
Consequently we get
\begin{eqnarray*}
\lefteqn{\langle N,\hat\nabla_{X_u}X_u-\hat\nabla_{X_v}X_v-2i\hat\nabla_{X_u}X_v\rangle} \\
& = & \frac1{1+|g|^2}(4g\bar F_{zz}+4\bar gF_{zz}+2ig\eta\bar F_z-2i\bar g\eta F_z+2(1-|g|^2)\eta_z).
\end{eqnarray*}

We first consider a domain where $g\neq 0$. Using \eqref{diffF} we obtain
$$\langle N,\hat\nabla_{X_u}X_u-\hat\nabla_{X_v}X_v-2i\hat\nabla_{X_u}X_v\rangle
=i\eta^2+2\frac{\bar g_z\eta}{\bar g}.$$ Finally, using \eqref{eta} we get
$$\left(\langle N,\hat\nabla_{X_u}X_u-\hat\nabla_{X_v}X_v-2i\hat\nabla_{X_u}X_v\rangle-i\eta^2\right)\rmd z^2
=4iQ.$$

We now deal with the points where $g=0$. At these points we have 
$\langle N,\hat\nabla_{X_u}X_u-\hat\nabla_{X_v}X_v-2i\hat\nabla_{X_u}X_v\rangle
=2\eta_z$, $\eta_z=8i\frac{g_z\bar g_z}{(1-|g|^2)^2}$ by \eqref{eta} and $\eta=0$. Hence the result also holds.
\end{proof}

\begin{rem}
We have the following formula for the Hopf differential in terms of the Weierstrass data of the surface:
$$Q=-\frac i2\frac{\bar g_z\eta}{\bar g}\rmd z^2.$$
In the case of minimal surfaces in $\R^3$, if $g$ denotes the Gauss map (meromorphic) of the surface and $\eta=2(x_3)_z$, where $x_3$ is the third coordinate of the surface, then the Hopf differential of the surface is given by $$Q=\frac12\frac{g'\eta}g\rmd z^2.$$  
\end{rem}

\section{The action of isometries} \label{isometries}

Let $\Sigma$ be a simply-connected Riemann surface. Let $X:\Sigma\to\nil$ be a nowhere vertical conformal minimal immersion and $g:\Sigma\to\h^2$ its Gauss map (we assume that the normal is upward); it is harmonic and nowhere antiholomorphic.

Let $\Phi$ be an isometry of $\nil$ preserving the orientation of the fibers. If $\Phi$ is a translation, then the Gauss map of $\Phi\circ X$ is also $g$, since the canonical frame $(E_1,E_2,E_3)$ is invariant by translations. From this we deduce that, if we remove the initial condition in theorem \ref{weierstrass}, then the immersion with a prescribed Gauss map is unique up to translations of $\nil$. If $\Phi$ is a rotation about a fiber, then the Gauss map of $\Phi\circ X$ is $\rho\circ g$, where $\rho$ is the rotation about $0\in\h^2$ of the same angle. Thus the isometries of $\nil$ only induce the rotations about $0$ in $\h^2$ (we do not take into consideration the isometries of $\nil$ that do not preserve the orientation of the fibers since we only consider surfaces whose normal vector is upward).

Conversely, let $T$ be a positive (i.e., orientation-preserving) isometry of $\h^2$. Then $T\circ g$ is harmonic and nowhere antiholomorphic; hence it is the Gauss map of a nowhere vertical conformal minimal immersion $\tilde X:\Sigma\to \nil$. 

If $T$ is a rotation about $0\in\h^3$, then $\tilde X$ is simply $\Phi\circ X$ where $\Phi$ is a rotation about a fiber in $\nil$. On the contrary, 
if $T$ is not a rotation about $0$, then in general $\tilde X$ is a non-trivial deformation of $X$. Example \ref{transl} illustrates this fact.

Consequently, a nowhere vertical conformal minimal immersion has a natural $2$-parameter family of deformations obtained by applying an isometry of $\h^2$ to the Gauss map. These deformations are in general not isometric. However, if $g$ has a positive symmetry, i.e. if there exists a positive isometry $T$ of $\h^2$ and a positive conformal diffeomorphism $\psi$ of $\Sigma$ such that $T\circ g=g\circ\psi$, then $g$ and $T\circ g$ define the same surface up to a translation in $\nil$ and up to a reparametrization.


If $T$ is a negative (i.e., orientation-reversing) isometry of $\h^2$, then $T\circ g$ is still harmonic but not necessarily nowhere antiholomorphic. Even if it is nowhere antiholomorphic, the corresponding minimal surface can be completely different from the one we started with. In particular, completeness may not be preserved. This fact is illustrated by examples \ref{semitrough} and \ref{semitrough2}.

\section{Complete minimal surfaces} \label{sectioncomplete}

The object of this section is to study the image of the Gauss map of complete minimal surfaces in $\nil$. For this purpose, we recall a few facts on harmonic maps and CMC surfaces in Minkowski space $\LL^3$.

Let $\Sigma$ be a simply connected Riemann surface, $X:\Sigma\to\nil$ a nowhere vertical minimal immersion and $g:\Sigma\to\h^2$ its Gauss map. The map $g$ is harmonic and nowhere antiholomorphic; hence it is the Gauss map of a spacelike CMC $\frac12$ immersion $\hat X:\Sigma\to\LL^3$, and the metric induced by $\hat X$ on $\Sigma$ is
$$\rmd\hat s^2=\frac{16}{(1-|g|^2)^2}|g_z|^2|\rmd z|^2$$
(see for example \cite{akni} and \cite{wan}; in \cite{akni}, $g_{\bar z}$ appears instead of $g_z$ because of a different choice of complex structure for $\h^2$; in \cite{wan}, the mean curvature is normalized to $1$ instead of $\frac12$). By \eqref{metric} we observe that the metric induced by $X$ on $\Sigma$ is
\begin{equation} \label{twometrics}
\rmd s^2=\left(\frac{1+|g|^2}{1-|g|^2}\right)^2\rmd\hat s^2.
\end{equation}
Hence, if $\rmd\hat s^2$ is complete, then $\rmd s^2$ is also complete.

Using these facts and some results about harmonic maps and CMC surfaces in $\LL^3$, we will deduce some properties of the image of the Gauss map of complete nowhere vertical minimal surfaces in $\nil$.

Let $\partial_{\infty}\h^2$ be the asymptotic boundary of the hyperbolic plane; we will identify it with $\s^1$. We will consider the intersection of $\partial_{\infty}\h^2$ and of the closure of the image of the Gauss map. This intersection is the set of all the asymptotic horizontal directions of the normal of the surface. 

\begin{thm} \label{gaussianimage}
1. The Gauss map of a complete nowhere vertical minimal surface $\cS$ in $\nil$ is not bounded in $\h^2$.

2. Let $L\subset\s^1=\partial_\infty\h^2$ be a closed subset containing at least $3$ points. Then there exists a complete nowhere vertical minimal surface in $\nil$ whose Gauss map is a diffeomorphism onto the interior of the convex hull of $L$ in $\h^2$.

3. If $L$ moreover has non-empty interior, then there exists a minimal surface as in 2. that is moreover of hyperbolic conformal type.

4. Let $m\geqslant 3$ and $\cP$ be an ideal polygon in $\h^2$ with $m$ vertices. Then there exists a $(m-3)$-parameter family of non-congruent complete nowhere vertical minimal surfaces in $\nil$ of parabolic conformal type and whose Gauss maps are diffeomorphisms onto the interior of $\cP$.
\end{thm}

\begin{proof}
1. Let $X:\Sigma\to\nil$ be a conformal immersion such that $X(\Sigma)=\cS$ and $g:\Sigma\to\h^2$ its Gauss map; without loss of generality we can assume that $\Sigma$ is simply-connected. Let $\hat X:\Sigma\to\LL^3$ be a CMC $\frac12$ immersion with Gauss map $g$; we use the same notations as above.

Assume that $g$ is bounded, i.e., that there is a constant $c\in(0,1)$ such that $|g|\leqslant c$. Then the quantity $\left(\frac{1+|g|^2}{1-|g|^2}\right)^2$ is bounded from above by a constant, and since $\rmd s^2$ is complete, so is $\rmd\hat s^2$ by \eqref{twometrics}. By a theorem of Y. L. Xin (\cite{xin}), the only complete CMC surface in $\LL^3$ with bounded Gauss map is the plane, which has CMC $0$. This gives a contradiction since $\hat X$ has CMC $\frac12$.

2. By theorem 4.7 in \cite{choitreibergs}, there exists a complete spacelike CMC surface given by a conformal immersion $\hat X:\Sigma\to\LL^3$ (where $\Sigma$ is simply-connected) such that the closure of the image of its Gauss map $g:\Sigma\to\h^2$ is the convex hull of $L$. Since $L$ contains at least three points, the map $g$ is a diffeomorphism by theorem 4.8 in \cite{choitreibergs}. Up to a change of orientation of $\Sigma$, we can assume that $g$ preserves orientation and so is nowhere antiholomorphic. Then it is the Gauss map of a nowhere vertical conformal minimal immersion $X:\Sigma\to\nil$. And $X(\Sigma)$ is complete since $\hat X(\Sigma)$ is complete. Then the surface $X(\Sigma)$ has the required properties.

3. We proceed as in 2. using theorem 9.2 in \cite{choitreibergs}. 

4. We proceed as in 2. using \cite{httw}. The fact that the obtained surfaces are non-congruent is a consequence of the fact that their Gauss maps are nontrivially distinct. 
\end{proof}

\begin{rem}
I. Fern\'andez and P. Mira proved that any holomorphic quadratic differential on $\C$ or $\D$ is the Abresch-Rosenberg differential of a complete minimal surface in $\nil$ (\cite{fernandezmira}).
\end{rem}

We end with some open questions.
It would be interesting to study the complete minimal surfaces with these prescribed Gauss maps. We know that they are local graphs over parts of $\R^2$. Hence it is a natural question to determine if they are proper, if they are graphs. It also seems interesting to try to solve the Bernstein problem in $\nil$ using the harmonic Gauss map, i.e., to classify all minimal entire graphs (i.e., over the entire $\R^2$) in terms of their Gauss map (example \ref{semitrough} is a new example of such a graph). Another problem is to determine if some complete nowhere vertical minimal surfaces in $\nil$ come from non-complete CMC surfaces in $\LL^3$.

\section{examples} \label{sectionexamples}

In this section we give the Gauss map of some classical examples and we construct new examples, in a particular a new example of an entire minimal graph over $\R^2$ (example \ref{semitrough}). Vertical planes do not appear, since their normal is horizontal (and constant).

\begin{exa}[the minimal hemisphere] \rm \label{hemisphere}
Let $\Sigma=\D$ and $g(z)=iz$. Then $g_z=i$, and using the formulas in theorem \ref{weierstrass} we get (up to a translation)
$$F=\frac{4z}{1-|z|^2},\quad h=0.$$ Hence we obtain the surface of equation $x_3=0$. It is invariant by rotations about the $x_3$-axis (but by no translation nor other rotation). We call it the \emph{minimal hemisphere} of pole $(0,0,0)$. We chose this terminology, since this surface is the limit of the half CMC spheres of South (or North) pole $(0,0,0)$ when the mean curvature tends to $0$. The Hopf differential of $g$ is $Q=0$.
\end{exa}

\begin{exa}[translation-invariant examples] \rm \label{transl}
Assume that the image of $g$ is a geodesic of $\h^2$ and that $g$ is nowhere antiholomorphic. Then the Hopf differential of $g$ cannot vanish (otherwise we would have $\bar g_z=0$ at some point, and since $g$ only depends on one real parameter this would lead to $g_z=0$ at this point, which is excluded); hence, up to a conformal change of parameter, we can choose $\Sigma=\C$, $Q=-\frac14\rmd z^2$ and
$$g(u+iv)=\frac{a(1-\cosh v)+\bar c\sinh v}{c(1-\cosh v)+\bar a\sinh v}$$ with $|a|^2-|c|^2=1$. Then we have
$$g_z=-\frac i2\frac{1-\cosh v}{(c(1-\cosh v)+\bar a\sinh v)^2}.$$ Up to a rotation about $0\in\h^2$ (see the discussion in section \ref{isometries}), we can assume that the geodesic intersects the axis $\{\re\zeta=0\}$ orthogonally and is oriented from the right to the left; moreover we can assume that this intersection happens when $v=0$. From this we get $g(0)\in i\R$ and $g_z(0)=-\frac i2g_v(0)\in i\R_+$, which gives $a\in\R$ and $c\in i\R$. Hence, up to a multiplication by $-1$ of $a$ and $c$, we can set $a=\cosh\theta$ and $c=i\sinh\theta$ for some $\theta\in\R$. From the formulas of theorem \ref{weierstrass} we get
$$F_z=\frac12(\cosh(2\theta)+\cosh v+i\sinh(2\theta)\sinh v),$$
$$F_{\bar z}=\frac12(\cosh(2\theta)-\cosh v-i\sinh(2\theta)\sinh v),$$
and so (up to a translation)
$$F=\cosh(2\theta)u-\sinh(2\theta)\cosh v+i\sinh v.$$ Then we compute that
$$h_z=\frac14\cosh(2\theta)\sinh v-\frac i4\sinh(2\theta)-\frac i4\cosh(2\theta)u\cosh v,$$
and finally we obtain (up to a translation)
$$h=\frac12\cosh(2\theta)u\sinh v+\frac12\sinh(2\theta)v.$$ This surface is an entire graph over $\R^2$, given by
$$x_3=\frac12x_1x_2+\frac12\sinh(2\theta)\left(\argsinh(x_2)+x_2\sqrt{1+x_2^2}\right).$$
It is invariant by the one-parameter family of translations along the $x_1$-axis. It is described in \cite{mercuri}.
\end{exa}

\begin{exa}[helicoids] \rm \label{helicoid}
Let $a\in\R\setminus\{0\}$ and $\varphi$ a real-valued function defined on an open interval $I$ satisfying $$(\varphi')^2-\varphi^2=a(1-\varphi^2)^2$$
and constant on no interval. We have $\varphi"-\varphi=-2a\varphi(1-\varphi^2)$, and so
$$(1-\varphi^2)\varphi"+2\varphi(\varphi')^2-\varphi(1+\varphi^2)=0.$$

Assume first that $a>0$. We choose $\varphi$ such that $\varphi(0)=0$ and $\varphi'>0$, and we restrict $\varphi$ to an interval $(-v_0,v_0)$ such that $\varphi(-v_0,v_0)=(-1,1)$. (If $a=\frac12$, then $\varphi(x)=\tan\left(\frac x2\right)$.)
Let $\Sigma=\R\times(-v_0,v_0)$ and $g(u+iv)=e^{iu}\varphi(v)$; then $g$ is harmonic.
We have $$g_z=\frac i2e^{iu}(\varphi-\varphi'),$$ and so
$$F_z=2e^{iu}\frac{\varphi-\varphi'}{(1-\varphi^2)^2},\quad
F_{\bar z}=-2e^{iu}\varphi^2\frac{\varphi-\varphi'}{(1-\varphi^2)^2}.$$
Up to a translation we get $$F=-2ie^{iu}\frac{\varphi(v)-\varphi'(v)}{1-\varphi(v)^2}.$$
Then we compute that $$h_z=\frac{\varphi'(v)^2-\varphi(v)^2}{(1-\varphi(v)^2)^2}=a,$$ which gives, up to a translation,
$$h=au.$$ We have $|F|\to+\infty$ when $v\to-v_0$, and $|F|\to 0$ when $v\to v_0$. Thus the closure of the surface contains the $x_3$-axis; it is one half of a right-helicoid. We obtain the complete right-helicoid by rotation of angle $\pi$ about this axis. Observe that the curve $\{v=0\}$ on the surface is horizontal and the helicoid is horizontal along this curve (since $g(u,0)=0$). The Hopf differential of $g$ is $Q=-a\rmd z^2$.

Assume now that $a<0$. We choose $\varphi$ such that $\varphi>0$ and $\varphi$ has a minimum at $v=0$, and we restrict $\varphi$ to an interval $(-v_0,v_0)$ such that $\varphi(-v_0,v_0)=(\varphi(0),1)$. Let $\Sigma=\R\times(-v_0,v_0)$ and $g(u+iv)=e^{iu}\varphi(v)$. As above, $g$ is harmonic and the same formulas for $F$ and $h$ hold.
In the same way, we have $|F|\to+\infty$ when $v\to-v_0$ and $|F|\to 0$ when $v\to v_0$. Thus the closure of the surface contains the $x_3$-axis; it is one half of a left-helicoid. We obtain the complete left-helicoid by rotation of angle $\pi$ about this axis. The Hopf differential of $g$ is $Q=-a\rmd z^2$.
\end{exa}

\begin{exa}[an entire minimal graph] \rm \label{semitrough}
Let $\Sigma=\{z\in\C;\re z>0\}$ and $$g(u+iv)=\frac{-1+i\cosh(2u)\sinh(2v)}{\sinh(2u)+\cosh(2u)\cosh(2v)}.$$ This is the Gauss map of a complete CMC surface in $\LL^3$ called semitrough (see \cite{choitreibergs} and \cite{httw} page 98). The image of $g$ is $\h^2\cap\{\re z<0\}$. We get 
$$g_z=\frac{(\cosh(2u)+1)(\cosh(2u)+\sinh(2u)\cosh(2v)-i\sinh(2v))}{(\sinh(2u)+\cosh(2u)\cosh(2v))^2},$$ and so
$$F_z=-\frac{i\cosh(2u)+\sinh(2v)+i\cosh(2v)\sinh(2u)}{\cosh(2u)-1},$$
$$F_{\bar z}=-\frac{i\cosh(2u)+\sinh(2v)-i\cosh(2v)\sinh(2u)}{\cosh(2u)-1}.$$
Up to a translation we obtain
$$F=2\sinh v\cosh{v}\coth{u}+i\coth u-2iu,$$
and then 
\begin{eqnarray*}
h_z & = & i\cosh^2v-\frac i2-u\frac{\sinh v\cosh v}{\sinh^2u}+iu\coth u \\
& & -2iu\cosh^2v\coth u+\sinh v\cosh v\coth u,
\end{eqnarray*}
and finally, up to a translation,
$$h=\sinh v\cosh v(2u\coth u-1).$$ The map $F$ is a diffeomorphism from $\Sigma$ onto $\C$. Hence we obtain a new example of an entire graph over $\R^2$. Let $x_1=\re F$, $x_2=\im F$ and $x_3=h$ be the coordinates of the surface. The map $u\mapsto x_2$ is a diffeomorphism from $(0,+\infty)$ onto $\R$, and the surface is defined by an equation of the form $$x_3=x_1f(x_2),$$ where $f$ is a function. This surface is foliated by Euclidean straight lines (which, in general, are not geodesics of $\nil$). When $x_2\to-\infty$, i.e., when $u\to+\infty$, we have $f(x_2)\sim-\frac12x_2$; when $x_2\to+\infty$, i.e., when $u\to0$, we have $f(x_2)\sim\frac1{2x_2}$. Hence on the one side the surface is asymptotic to the surface of equation $x_3=-\frac12x_1x_2$, which is invariant by translations along the $x_2$-axis (this is the image of example \ref{transl} with $\theta=0$ by a rotation); on the other side, the surface is asymptotic to the surface of equation $x_3=0$, which is a minimal hemisphere (example \ref{hemisphere}). This surface is complete, since the corresponding CMC $\frac12$ surface in $\LL^3$ is complete. The Hopf differential of $g$ is $Q=-\rmd z^2$.
\end{exa}

\begin{ex} \label{semitrough2}
Let $\Sigma=\{z\in\C;\re z>0\}$ and $$g(u+iv)=\frac{-1-i\cosh(2u)\sinh(2v)}{\sinh(2u)+\cosh(2u)\cosh(2v)}.$$ This is the complex conjugate of the Gauss map of example \ref{semitrough}. We get $$g_z=-\frac{(\cosh(2u)-1)(\cosh(2u)+\sinh(2u)\cosh(2v)+i\sinh(2v))}{(\sinh(2u)+\cosh(2u)\cosh(2v))^2},$$ and so
$$F_z=\frac{i\cosh(2u)-\sinh(2v)+i\cosh(2v)\sinh(2u)}{\cosh(2u)+1},$$
$$F_{\bar z}=\frac{i\cosh(2u)-\sinh(2v)-i\cosh(2v)\sinh(2u)}{\cosh(2u)+1}.$$
Up to a translation we obtain
$$F=-2\sinh v\cosh{v}\tanh{u}-i\tanh u+2iu,$$
and then 
\begin{eqnarray*}
h_z & = & i\cosh^2v-\frac i2+u\frac{\sinh v\cosh v}{\cosh^2u}+iu\tanh u \\
& & -2iu\cosh^2v\tanh u+\sinh v\cosh v\tanh u,
\end{eqnarray*}
and finally, up to a translation,
$$h=\sinh v\cosh v(2u\tanh u-1).$$ The map $F$ is a diffeomorphism from $\Sigma$ onto $\R\times(0,+\infty)$. Hence the surface is a graph over a half-plane. Let $x_1=\re F$, $x_2=\im F$ and $x_3=h$ be the coordinates of the surface. The map $u\mapsto x_2$ is a diffeomorphism from $(0,+\infty)$ onto itself, and the surface is defined by an equation of the form $$x_3=x_1f(x_2),$$ where $f$ is a function. This surface is foliated by Euclidean straight lines (as example \ref{semitrough}). When $x_2\to+\infty$, i.e., when $u\to+\infty$, we have $f(x_2)\sim-\frac12x_2$; when $x_2\to0$, i.e., when $u\to0$, we have $f(x_2)\sim\frac1{2x_2}$. Hence on the one side the surface is asymptotic to the surface of equation $x_3=-\frac12x_1x_2$, and on the other side the surface is not complete and its closure contains the $x_3$-axis. The surface can be completed by rotation of angle $\pi$ about the $x_3$-axis. The Hopf differential of $g$ is $Q=-\rmd z^2$.
\end{ex}

\section{Final remark}

After the first version of this paper ({arXiv:math/0606299v1}) was written, in June 2006, some of the questions addressed at the end of section \ref{sectioncomplete} were answered:
\begin{itemize}
\item in \cite{fmtams}, I. Fern\'andez and P. Mira solved the Bernstein problem in $\nil$, i.e., they classified all entire minimal graphs in $\nil$ (in terms of their Abresch-Rosenberg differential),
\item in \cite{dh}, the author and L. Hauswirth proved that every complete nowhere vertical minimal surface in $\nil$ is necessarily an entire graph; moreover, this implies using \cite{fmtams} that every complete nowhere vertical minimal surface in $\nil$ comes from a complete CMC surface in $\LL^3$.
\end{itemize}

\bibliographystyle{alpha}
\bibliography{gauss}

\end{document}